\documentclass{amsart}
 \usepackage{amscd}
                     \usepackage{amssymb}
                     \usepackage{amsmath}
\usepackage{amsthm}
\usepackage{url}

 \newcommand{\be}{\begin{equation}}
       \newcommand{\ee}{\end{equation}}
       \newcommand{\ba}{\begin{eqnarray}}
        \newcommand{\ea}{\end{eqnarray}}
 \newcommand{\ban}{\begin{eqnarray*}}
 \newcommand{\ean}{\end{eqnarray*}}

 \newcommand{\lp}{\langle}
 \newcommand{\rp}{\rangle}
 \newcommand{\ra}{\rightarrow}

 \newcommand{\sect}[1]{\section{#1} \setcounter{equation}{0}}

 \newcommand{\Ric}{\mathrm{Ric}}

 \newcommand{\Hess}{\mathrm{Hess}}
  \newcommand{\tr}{\mbox{tr}}
 
 \DeclareMathOperator{\divg}{div}

 \newtheorem{theo}{Theorem}[section]

\begin{document}
 \newtheorem{defn}[theo]{Definition}
 \newtheorem{ques}[theo]{Question}
 \newtheorem{lem}[theo]{Lemma}
 \newtheorem{lemma}[theo]{Lemma}
 \newtheorem{prop}[theo]{Proposition}
 \newtheorem{coro}[theo]{Corollary}
 \newtheorem{ex}[theo]{Example}
 \newtheorem{note}[theo]{Note}
 \newtheorem{conj}[theo]{Conjecture}
 \newtheorem{Ex}[theo]{Example}
 \newtheorem{Remark}[theo]{Remark}

\title{ Rigidity of Quasi-Einstein Metrics}
\author{Jeffrey Case}
\email{casej@math.ucsb.edu}
 \author{Yu-Jen Shu}
 \email{yjshu@math.ucsb.edu}
 \author{Guofang Wei}
 \email{wei@math.ucsb.edu}
 \address{Department of Mathematics, UCSB, Santa Barbara, CA 93106}
  \thanks{Partially supported by NSF grant DMS-0505733}
\date{}
\maketitle
\begin{abstract}
We call a metric quasi-Einstein if the $m$-Bakry-Emery Ricci tensor is a constant multiple of the metric tensor. This is a generalization of  Einstein metrics, which contains gradient Ricci solitons and is also closely related to the construction of the warped product Einstein metrics. 
We  study properties of  quasi-Einstein metrics and prove several rigidity  results. We also give a splitting theorem for some K\"ahler quasi-Einstein metrics. \end{abstract}

\sect{Introduction}

Einstein metrics and their generalizations are important both in mathematics and physics. A particular example is from the study of smooth metric measure spaces. 
Recall a smooth metric measure space is a  triple $(M^n,g, e^{-f} dvol_g)$, where $M$ is a complete
$n$-dimensional Riemannian manifold with metric $g$, $f$ is a smooth
real valued function on $M$, and $dvol_g$ is the Riemannian volume
density on $M$.  A natural extension of the Ricci tensor to smooth metric measure
spaces is the $m$-Bakry-Emery Ricci tensor
\be \Ric_f^m = \Ric +\Hess f  -\frac 1m df \otimes df  \ \ \  \mbox{for} \
0< m \le \infty. \label{Ric-m-f} \ee
When $f$ is constant, this is the usual Ricci tensor. We call a triple $(M, g, f)$ (a Riemannian manifold  $(M,g)$ with a function $f$ on $M$) (m-)quasi-Einstein if it
 satisfies the  equation
\begin{equation} \label{sol} \Ric_f^m  =\Ric +\Hess f  -\frac 1m df \otimes df = \lambda g \end{equation}
for some $\lambda \in\mathbb{R}$.  This equation is especially interesting in that when $m = \infty$ it is exactly the gradient Ricci soliton equation; 
when $m$ is a positive integer, it corresponds to warped product Einstein metrics (see Section~\ref{warped} for detail); when $f$ is constant,  it 
gives the Einstein equation. We call a quasi-Einstein metric trivial when $f$ is constant (the rigid case).

Many geometric and topological properties of manifolds with Ricci curvature bounded below can be extended to manifolds with $m$-Bakry-Emery Ricci 
tensor bounded from below when $m$ is finite or $m$ is infinite and $f$ is bounded, see the survey article \cite{Wei-Wylie2007} and the references 
there for details. 
 
Quasi-Einstein metrics for  finite $m$ and for $m= \infty$ share some common properties. It is well-known now that compact solitons with $\lambda \le 0$ are trivial \cite{Ivey1993}.
The same result is proven in \cite{KK} for quasi-Einstein metrics on compact manifolds with finite $m$.  Compact shrinking Ricci solitons have positive scalar 
curvature \cite{Ivey1993, Eminenti-Nave-Mantegazza2006}. Here we show 
\begin{prop}
A quasi-Einstein metric with  $1 \le m <  \infty$ and  $\lambda >0$ has positive scalar curvature.
\end{prop}

In dimension 2 and 3 compact Ricci solitons are trivial \cite{Hamilton1988, Ivey1993}. More generally compact shrinking  Ricci solitons with zero Weyl tensor are trivial \cite{Eminenti-Nave-Mantegazza2006,PW2, Ni-Wallach07}. We prove a similar result in dimension 2 for $m$ finite.
\begin{theo}  \label{dim2}
All 2-dimensional quasi-Einstein metrics on compact  manifolds are trivial.
\end{theo}

In fact, from the correspondence with warped product metrics, 2-dimensional quasi-Einstein metrics (with finite $m$) can be classified, see \cite[Theorem 9.119]{Besse}. (The  proof was not published though.)

 In Section~\ref{formula}  we also extend several properties for Ricci solitons ($m= \infty$) to quasi-Einstein metrics (general $m$). 
  
On the other hand  we show K\"ahler quasi-Einstein metrics behave very differently when $m$ is finite and $m$ is infinite.  

\begin{theo}  \label{Kahler}
Let $(M^n, g)$ be an n-dimensional complete simply-connected Riemannian manifold with a 
K\"ahler quasi-Einstein metric for finite $m$. Then $M=M_1\times M_2$ is a
Riemannian product, and $f$ can be considered as a function of $M_2$,
where $M_1$ is an (n-2)-dimensional Einstein manifold with Einstein
constant $\lambda$, and $M_2$ is a 2-dimensional quasi-Einstein
manifold.
\end{theo}

Combine this with Theorem~\ref{dim2} we immediately get 
\begin{coro}  \label{coro}
There are no nontrivial K\"ahler quasi-Einstein metrics with  finite $m$ on compact manifolds.
\end{coro}

Note that  all  known
 examples of (nontrivial) compact shrinking soliton are K\"ahler, see the survey article \cite{Cao2006}.  

Ricci solitons play a very important role in the theory of Ricci flow and  are extensively studied recently. 
Warped product Einstein metrics have considerable interest in physics and
many Einstein metrics are constructed in this form, especially on noncompact
manifolds \cite{Besse}. It was asked in \cite{Besse} whether one could find Einstein metrics with nonconstant warping function on compact manifolds. From Corollary~\ref{coro} only 
non-K\"ahler ones are possible.  Indeed in \cite{Lu-Page-Pope2004}  warped product Einstein metrics are constructed on a class of $S^2$ bundles over K\"ahler-Einstein 
bases warped with $S^m$ for $m \ge2$, giving compact nontrivial quasi-Einstein metrics for positive integers $n\ge 4$,  $m \ge 2$. When $m=1$, there are no nontrivial  quasi-Einstein metrics on compact manifolds, see Remark~\ref{m=1} and Proposition~\ref{R-const}.  When $n=3, m \ge 2$ it remains open. 

In Section 4 we also give a characterization of  quasi-Einstein metrics with finite m which are Einstein at the same time.

\section{Warped Product Einstein Metrics}  \label{warped}

In this section we show that when $m$ is a positive integer the quasi-Einstein metrics \eqref{sol} correspond to some warped product Einstein metrics, mainly due to the work of \cite{KK}.

Recall that given two Riemannian manifolds $(M^n,g_M)$,  $(F^m,g_F)$ and a positive smooth function $u$ on $M$, the warped product metric on $M\times F$  is defined by
\be
g = g_M + u^2 g_F.
\ee
We denote it as $M\times_u F$. Warped product is very useful in constructing various metrics.

When $0<m<\infty$, consider $u = e^{-\frac fm}$. Then we have
\ban
\nabla u & = & -\frac 1m \,  e^{-\frac fm} \nabla f, \\
\frac mu \, \Hess\, u & = & -\Hess f + \frac 1m df \otimes df. \ean
Therefore  \eqref{sol} can be rewritten as  \be  \label{sol-u} \Ric -
\frac mu \,\Hess\, u = \lambda g. \ee Hence we can use equation
\eqref{sol-u} to study \eqref{sol} when $m$ is finite and vice verse. Taking trace of
\eqref{sol-u} we have \be \Delta u = \frac um \left(R- \lambda
n\right). \ee Since $u >0$ this immediately gives the following
result which is similar to the $m =\infty$ (soliton) case.
\begin{prop}  \label{R-const}
A compact quasi-Einstein metric with constant scalar curvature is
trivial.
\end{prop}

In \cite{KK} it is shown that a Riemannian manifold (M,g) satisfies
\eqref{sol-u} if and only if the warped product metric $M\times_u
F^m$ is Einstein, where $F^m$ is an $m$-dimensional Einstein
manifold with Einstein constant $\mu$ satisfies
\[
\mu = u \Delta u +(m-1) |\nabla u|^2 + \lambda u^2.
\]
(In \cite{KK} it is only stated for compact Riemannian manifold, while
compactness is redundant. Also the Laplacian there and here have
different sign.) Therefore we have the following nice
characterization of the quasi-Einstein metrics as the base metrics of warped product
Einstein metrics.
\begin{theo}
 $(M,g)$ satisfies the quasi-Einstein equation  \eqref{sol} if and only if
 the warped product metric $M\times_{e^{-\frac fm}} F^m$ is Einstein, where $F^m$ is an $m$-dimensional Einstein manifold with Einstein constant $\mu$ satisfying
 \begin{equation} \label{mu} \mu e^{\frac{2}{m}f} = \lambda - \frac{1}{m}\left(\Delta f - |\nabla f|^2\right). \end{equation}
\end{theo}

\section{Formulas and Rigidity for Quasi-Einstein Metrics}  \label{formula}
In this section we generalize the calculations in \cite{PW} for Ricci solitons to the metrics satisfying the  quasi-Einstein equation (\ref{sol}).

Recall the following general formulas, see e.g. \cite[Lemma 2.1]{PW} for a proof.
\begin{lem}  For a function $f$ in a Riemannian manifold
\be
\label{e2} 2(\divg\Hess f)(\nabla f)  =  \frac{1}{2}\Delta|\nabla f|^2 - |\Hess f|^2 + \Ric(\nabla f,\nabla f) + \lp\nabla f,\nabla\Delta f\rp,  
\ee
\be
\label{e1} \divg\nabla\nabla f  = \Ric\nabla f + \nabla\Delta f.
\ee
\end{lem}
\vspace*{.2in}

The trace form of \eqref{sol}
\be  \label{sol-trace}
R + \Delta f - \frac 1m |\nabla f|^2 = \lambda n,
\ee
where $R$ is the scalar curvature, will be  used later.

Using these formulas and the contracted second Bianchi identity
\be \label{e3} \nabla R  = 2\divg \Ric, \ee
we can show the following formulas for quasi-Einstein metrics, which generalize some of the formulas in Section 2 of \cite{PW}.
\begin{lem} \label{formulas} If  $\Ric_f^m=\lambda g$,  then
\ba
\label{e4} \frac{1}{2}\Delta|\nabla f|^2  & = &  |\Hess f|^2 - \Ric(\nabla f,\nabla f) + \frac{2}{m}|\nabla f|^2\Delta f,   \\
\label{e5} \frac{1}{2}\nabla R & = & \frac{m-1}{m}\Ric(\nabla f) + \frac{1}{m}\left(R-(n-1)\lambda\right)\nabla f ,
\ea
\ba
\label{e8} \frac{1}{2}\Delta R - \frac{m+2}{2m}\nabla_{\nabla f}R & = & \frac{m-1}{m}\tr\left(\Ric\circ(\lambda I-\Ric)\right) - \frac{1}{m}(R-n\lambda)(R-(n-1)\lambda) \nonumber \\
\label{e9} & = & -\frac{m-1}{m}\left|\Ric-\frac{1}{n}Rg\right|^2 - \frac{m+n-1}{mn}(R-n\lambda)(R-\frac{n(n-1)}{m+n-1}\lambda).
\ea
\end{lem}
\begin{proof} From \eqref{e2} we have
\be \label{Delta} \frac{1}{2}\Delta|\nabla f|^2 =
2(\divg\Hess f)(\nabla f) + |\Hess f|^2 - \Ric(\nabla f,\nabla f) - \lp\nabla f,\nabla\Delta f\rp.
\ee
By taking the divergence of \eqref{sol}, we have
\be
 \divg\Ric + \divg\Hess f - \frac{1}{m}\Delta f\;df - \frac{1}{m}(\nabla_{\nabla f}\nabla f)^\ast = 0,  \label{divRic} \ee
where $X^\ast$ is the dual 1-form of  the vector field $X$. Using (\ref{e3})  we get
\be  \label{divHess}
2\divg\Hess f (\nabla f)  = -  \lp \nabla R, \nabla f \rp  +  \frac{2}{m}\Delta f\; |\nabla f|^2 + \frac{2}{m}\Hess f(\nabla f,\nabla f).
\ee
 Now taking the covariant derivative of (\ref{sol-trace}) yields
\begin{equation}\label{grad} \nabla R + \nabla\Delta f - \frac{1}{m}\nabla|\nabla f|^2 = 0. \end{equation}
Plug this into \eqref{divHess} and then plug \eqref{divHess} into (\ref{Delta}) we get
\ban
 \frac{1}{2}\Delta|\nabla f|^2  & = &  \lp \nabla\Delta f - \frac{1}{m}\nabla|\nabla f|^2, \nabla f \rp  +  \frac{2}{m}\Delta f\; |\nabla f|^2 + \frac{2}{m}\Hess f(\nabla f,\nabla f) \\
 & &  + |\Hess f|^2 - \Ric(\nabla f,\nabla f) - \lp\nabla f,\nabla\Delta f\rp \\
 & =& |\Hess f|^2 - \Ric(\nabla f,\nabla f) + \frac{2}{m}|\nabla f|^2\Delta f,
 \ean
which is  \eqref{e4}.

For  \eqref{e5},  using \eqref{e3}, \eqref{divRic}, \eqref{e1}, and \eqref{grad}, we get
\begin{align*}
\nabla R & = 2\divg\Ric \\
& = -2\divg\Hess f + \frac{2}{m}\Delta f\;\nabla f + \frac{2}{m}\nabla_{\nabla f}\nabla f \\
& = -2\Ric(\nabla f) - 2\nabla\Delta f + \frac{2}{m}\Delta f\;\nabla f + \frac{2}{m}\nabla_{\nabla f}\nabla f \\
& = -2\Ric(\nabla f)+2\nabla R - \frac{2}{m}\nabla|\nabla f|^2 + \frac{2}{m}\Delta f\;\nabla f + \frac{2}{m}\nabla_{\nabla f}\nabla f.
\end{align*}
Solving for $\nabla R$ and noting that $\nabla|\nabla f|^2 = 2\nabla_{\nabla f}\nabla f$, we get
\[ \nabla R = 2\Ric(\nabla f)- \frac{2}{m}\Delta f\;\nabla f + \frac{2}{m}\nabla_{\nabla f}\nabla f. \]
From \eqref{sol}, we have that
\[ \nabla_{\nabla f}\nabla f = \left(\lambda + \frac{1}{m}|\nabla f|^2\right)\nabla f - \Ric(\nabla f), \]
so by using this substitution and \eqref{sol-trace} for $\Delta f$, we arrive at \eqref{e5},
\begin{align*}
\frac{1}{2}\nabla R & = \frac{m-1}{m}\Ric(\nabla f) + \frac{1}{m}\left(R-(n-1)\lambda\right)\nabla f.
\end{align*}

Taking the divergent of the above equation we have
\be  \label{div}
\frac{1}{2}\Delta R  = \frac{m-1}{m} \divg \left( \Ric(\nabla f) \right)+ \frac{1}{m}\divg \left (\left(R-(n-1)\lambda\right)\nabla f\right).
\ee
Now
\ba
\divg \left( \Ric(\nabla f) \right) & = &  \lp \divg \Ric, \nabla f\rp + \tr \left(\Ric \circ \Hess f \right) \nonumber  \\
& = & \lp \frac 12 \nabla R, \nabla f\rp + \tr \left(\Ric \circ \left( \frac 1m df \otimes df  +\lambda g -\Ric \right) \right) \nonumber \\
& = & \lp \frac 12 \nabla R, \nabla f\rp + \frac 1m \Ric (\nabla f, \nabla f) + \tr \left(\Ric \circ \left(\lambda g -\Ric \right) \right)  \nonumber \\
& = &  \lp \frac 12 \nabla R, \nabla f\rp + \frac{1}{m-1} \left( \lp \frac 12 \nabla R, \nabla f\rp - \frac 1m \left( R-(n-1) \lambda \right) |\nabla f|^2 \right)  \label{divRic-f}  \\
& &  + \tr \left(\Ric \circ \left(\lambda g -\Ric \right) \right), \nonumber
\ea
where the last equation comes from \eqref{e5}. Also
\be  \label{divf}
\divg \left( \left(R-(n-1)\lambda\right)\nabla f \right)  =  \left(R-(n-1)\lambda\right)\Delta f +\lp \nabla R, \nabla f \rp.
\ee
Plugging \eqref{divRic-f} and \eqref{divf} into \eqref{div} and using \eqref{sol-trace}
 we arrive at
\[ \frac{1}{2}\Delta R - \frac{m+2}{2m}\nabla_{\nabla f}R = \frac{m-1}{m}\tr\left(\Ric\circ(\lambda I-\Ric)\right) - \frac{1}{m}(R-n\lambda)(R-(n-1)\lambda) . \]
Let $\lambda_i$ be the eigenvalues of the Ricci tensor, we get
\begin{align*}
\tr\left(\Ric\circ(\lambda I-\Ric)\right) & = \sum \lambda_i(\lambda-\lambda_i) \\
& = -\left|\Ric-\frac{1}{n}Rg\right|^2 + R\left(\lambda-\frac{1}{n}R\right),
\end{align*}
which yields \eqref{e9}. \end{proof}

As in \cite{PW}, these formulas give important information about quasi-Einstein metrics. Combining the first equation \eqref{e4} in Lemma~\ref{formulas} with the maximal principle we have
\begin{prop}
If  a compact Riemannian manifold satisfying \eqref{sol} and
\[
 \Ric(\nabla f,\nabla f) \le  \frac{2}{m}|\nabla f|^2\Delta f
 \] then the function $f$ is constant so  it is Einstein.
 \end{prop}

 Equation \eqref{e5} gives
\begin{prop}
When $m \not= 1$, a quasi-Einstein metric has constant scalar curvature if and only if
\[
 \Ric(\nabla f) = - \frac{1}{m-1} \left( R - (n-1)\lambda \right) \nabla f.
 \]
 \end{prop}

\begin{Remark}  \label{m=1} When $m=1$, the constant $\mu$ in \eqref{mu} is zero, combining with \eqref{sol-trace}, we get $R = (n-1)\lambda$. The scalar curvature is always constant.
\end{Remark}
 
 Equation \eqref{e9} gives the following results.
  \begin{prop}
If  a Riemannian manifold $M$ satisfies \eqref{sol} with $m \ge 1$ and

a) $\lambda >0$ and $M$  is compact then the scalar curvature is bounded below by
\be  \label{R-lower-b}
R \ge  \frac{n(n-1)}{m+n-1}\lambda.
\ee
Equality holds if and only if $m=1$.

b) $\lambda =0$,  the scalar curvature is constant and $m >1$, then $M$ is Ricci flat.

c) $\lambda < 0$ and the scalar curvature is constant, then \[
n \lambda \le R \le \frac{n(n-1)}{m+n-1}\lambda\]
and  when $m >1$, $R$ equals either of the extreme values iff $M$ is  Einstein.
 \end{prop}

\begin{Remark} When $m$ is finite, a manifold with quasi-Einstein metric and $\lambda >0$ is automatically compact \cite{Qian1997}. 
\end{Remark}

 \begin{Remark} Let $m= \infty$, we recover the well know result  \cite{Ivey1993} that compact shrinking Ricci soliton has positive scalar curvature, and some results in \cite{PW} about gradient Ricci solitons with constant scalar curvature.
 \end{Remark}
 
 \begin{proof}
a) Since $M$ is compact, applying the equation \eqref{e9} to a minimal point of $R$, we have
 \[
 - \frac{m+n-1}{mn}(R_{min}-n\lambda)(R_{min}-\frac{n(n-1)}{m+n-1}\lambda) \ge  \frac{m-1}{m}\left|\Ric-\frac{1}{n}Rg\right|^2 \ge0.
 \]
 So
 \[
\frac{n(n-1)}{m+n-1}\lambda \le  R_{min} \le n \lambda
\]
which gives \eqref{R-lower-b}.

b) c) Since $R$ is constant, from  \eqref{e9}
\[
 - \frac{m+n-1}{mn}(R-n\lambda)(R-\frac{n(n-1)}{m+n-1}\lambda) =  \frac{m-1}{m}\left|\Ric-\frac{1}{n}Rg\right|^2 \ge 0.\]
So if $\lambda = 0, \ m >1$, then $\Ric = \frac{1}{n}Rg$ and $R = 0$, thus it is Ricci flat. If $\lambda < 0$, $R \in [n\lambda,  \frac{n(n-1)}{m+n-1}\lambda]$.
 \end{proof}

\section{Two Dimensional Quasi-Einstein Metrics}

First we recall a characterization of warped product metrics found
in \cite{Cheeger-Colding} (see also \cite{PW2}).
\begin{theo}[Cheeger-Colding]   \label{Ch-Co-warp} A Riemannian manifold $(M^n,g)$ is a warped product $(a,b) \times_u N^{n-1}$
if and only if there is a nontrivial function $h$ such that $\Hess
\, h = k g$ for some function $k: M \ra \mathbb R$. ($u = h'$ up to a multiplicative constant) \end{theo}

From this we can give a characterization of quasi-Einstein metrics which are Einstein.
 \begin{prop}
 A complete finite $m$ quasi-Einstein metric $(M^n,g,u)$ is Einstein if and only if $u$ is constant or $M$ is diffeomorphic to $\mathbb R^n$ with the warped product structure $\mathbb R \times_{a^{-1}e^{ar}} N^{n-1}$, where $N^{n-1}$ is Ricci flat, $a$ is a constant (see below for its value).
 \end{prop}
 \begin{proof}
 If $g$ is Einstein, then $\Ric = \mu g$ for some constant $\mu$. From  \eqref{sol-u} we have
 \be  \label{Hess-u}
  \Hess\, u = \frac{\mu-\lambda}{m} u\, g \ \  \mbox{for some} \  u>0.  \ee
 If $M$ is compact, then $u$ (thus $f$) is constant. So if $u$ is not constant, then $M$ is noncompact and
$\lambda \le 0,  \mu \le 0$ and $\mu >  \lambda$.  So $\lambda <0,
\lambda < \mu \le 0$ and $u$ is a strictly convex function. Therefore $M^n$ is diffeomorphic to $\mathbb R^n$. By \eqref{Hess-u} and Theorem~\ref{Ch-Co-warp},  $u =c e^{\sqrt{\frac{\mu-\lambda}{m}} r}$, where $c$ is some constant. And $M$ is  $\mathbb R \times N^{n-1}$ with the warped product metric
\[
g = dr^2 + a^{-2}  e^{2a r}
g_0,
\]
where $a = \sqrt{\frac{\mu-\lambda}{m}}$. 
Since $g$ is Einstein we get  $ \mu = - (n-1) \frac{\mu-\lambda}{m} <0$ and $g_0$ is Ricci flat.
\end{proof}

\begin{Remark} The Taub-NUT metric \cite{Hawking1977} is a Ricci flat metric on $\mathbb R^4$ which is not flat.
\end{Remark}

Since a  2-dimensional Riemannian manifold satisfies $\Ric=\frac{R}{2}g$, we get  an  immediate  corollary of Theorem~\ref{Ch-Co-warp}.
\begin{coro}
A two dimensional quasi-Einstein metric \ref{sol-u} is a warped product metric.
\end{coro}

 Now we will prove Theorem~\ref{dim2}.
\begin{proof}
Since $M$ is compact, by \cite{KK}, we only need to prove the
theorem when $\lambda>0$. From \eqref{R-lower-b} we have \be
R\ge \frac{2}{m+1}\lambda. \ee So up to a cover we may assume $M$ is diffeomorphic to $S^2$.
 Since $M$ is 2-dimensional, we have
$\Ric=\frac{R}{2}g$.  Thus \eqref{e5} becomes \be \label{ee1} \nabla
R  = \frac{m+1}{m}\left(R-\frac{2}{m+1} \lambda \right)\nabla f. \ee

 Now let $u = e^{-\frac fm}$, then from \eqref{sol-u} $\Hess \, u= \frac um \left(\frac{R}{2}-\lambda\right)g$.
  In particular, $\nabla u$ is conformal.  By the Kazdan-Warner identity \cite{Kazdan-Warner1975}, we have
\[  \int_M\lp\nabla R,\nabla u\rp\;dV =0. \]
Thus \[
 -\frac 1m \int\lp\nabla R,\nabla f\rp e^{-\frac{f}{m}}\;dV =  0. \]
Using \eqref{ee1}, since $R\ge  \frac{2}{m+1}\lambda$,  we get  $\nabla f=\nabla R=0$.
\end{proof}

\section{K\"ahler Quasi-Einstein Metrics}
There are many nontrivial examples of shrinking K\"ahler-Ricci
solitons \cite{cao, wz}. In contrast, we will show here that
K\"ahler  quasi-Einstein metrics with finite $m$ are very rigid.
\newline

\noindent {\em Proof of Theorem~\ref{Kahler}} First, since on a K\"ahler manifold, the metric and
Ricci tensor  are both compatible with the complex structure $J$,
from \eqref{sol-u} we have  $\Hess \, u(JU,JV)= \Hess \, u(U, V)$
for all vector fields $U,V$. That implies \be \label{Jnabla}
J\nabla_X \nabla u=\nabla_{JX} \nabla u, \ee
 and the  $\phi$ defined by
$\phi(U,V)=\Hess \,u (JU, V)$ is an (1,1)-form. Note that $2 \Hess
\, u = L_{\nabla u}g$. By \eqref{Jnabla} $L_{\nabla u} JX = J
L_{\nabla u} X$, so $2 \phi = L_{\nabla u} \omega$, where $\omega$
is the K\"ahler form. Since $L_X d = d L_X$ and $\omega$ is closed we
have $\phi$ is closed. Furthermore, since the Ricci form is closed,
from \eqref{sol-u} we get that the (1,1)-form $\frac{\phi}{u}$ is
also closed, so $du\wedge \phi =0$. Now
$$ (du\wedge \phi)(U, V, W )=U u\cdot g(\nabla_{JV}\nabla u,
W)+V u\cdot g(\nabla_{JW}\nabla u, U)+W u\cdot g(\nabla_{JU}\nabla
u, V).$$
Let $U, V \perp\nabla u, W= \nabla u$, we have  $\Hess \, u(JU, V)=0$
for all $U, V \perp\nabla u$. Hence
\be  \label{nabla}
\nabla_X \nabla u \, \| J\nabla u \ \  \mbox{ for all} \ X\perp \nabla u.
\ee
Let $U= \nabla u, V = J\nabla u, W \perp \nabla u$ we get  $\nabla_{\nabla u}\nabla u \,  \| \nabla u$.

Now we consider the 2-dimensional distribution $T_1$ that is spanned
by $\nabla u$ and $J\nabla u$ at those points where $\nabla u$  is
nonzero.  We will show that $T_1$ = Span$\{\nabla u, J\nabla u\}$ is
invariant under parallel transport, i.e. if $\gamma$ is a path in $M$, and
$U$ is a parallel field along $\gamma$, then 
\begin{equation}\label{*}
\nabla_{\gamma '}\left(\frac{g(U,\nabla u)\nabla u}{|\nabla u|^2}+
\frac{g(U,J\nabla
u)J\nabla u}{|\nabla u|^2}\right)=0 .
\end{equation}
Since the covariant derivative is linear in $\gamma'$, we can prove
this in three cases:
\begin{enumerate}
\item  when
$\gamma '\bot\nabla u, J\nabla u$: so $\gamma '\bot\nabla u$ and $J\gamma '\bot\nabla u$, by \eqref{nabla}
$\nabla_{\gamma '} \nabla u=0$. Since the complex
structure $J$ is parallel,   $\nabla_{\gamma '} J\nabla
u=J\nabla_{\gamma '} \nabla u=0$.  By assumption $\nabla_{\gamma '}U= 0$, hence \eqref{*} follows.
\item when $\gamma '=\nabla u$: Using $\nabla_{\nabla u}\nabla u \|
\nabla u$, we have

\ban &&\nabla_{\gamma '}\left(\frac{g(U,\nabla u)\nabla u}{|\nabla
u|^2}+
\frac{g(U,J\nabla u)J\nabla u}{|\nabla u|^2}\right) \nonumber\\
&=&\gamma '\left(\frac{g(U,\nabla u)}{|\nabla u|^2}\right)\nabla u +
\frac{g(U,\nabla u)}{|\nabla u|^2}\nabla_{\gamma'}\nabla
u+\gamma'\left(\frac{g(U,J\nabla u)}{|\nabla u|^2}\right)J\nabla u +
\frac{g(U,J\nabla
u)}{|\nabla u|^2}\nabla_{\gamma'}J\nabla u\nonumber\\
&=&\frac{g(U,\nabla_{\nabla u}\nabla u)}{|\nabla u|^2}\nabla
u-\frac{2g(\nabla u,\nabla_{\nabla u}\nabla u)g(U,\nabla u)}{|\nabla
u|^4}\nabla u +  \frac{g(U,\nabla
u)}{|\nabla u|^2}\frac{g(\nabla u,\nabla_{\nabla u}\nabla u)}{|\nabla u|^2}\nabla u\nonumber \\
&+&\frac{g(U,\nabla_{\nabla u}J\nabla u)}{|\nabla u|^2}J\nabla
u-\frac{2g(\nabla u,\nabla_{\nabla u}\nabla u)g(U,J\nabla
u)}{|\nabla u|^4}J\nabla u + \frac{g(U,J\nabla
u)}{|\nabla u|^2}\frac{g(J\nabla u,\nabla_{\nabla u}J\nabla u)}{|\nabla u|^2}J\nabla u\nonumber \\
&=&\frac{\nabla u}{|\nabla u|^2}\left(\Hess \,  u(U, \nabla u)-\frac{\Hess \, u (\nabla u, \nabla u)}{|\nabla u|^2}g(U, \nabla u) \right)\nonumber\\
&+&\frac{J\nabla u}{|\nabla u|^2}\left(-\Hess \, u (JU, \nabla
u)+\frac{\Hess
u(\nabla u, \nabla u)}{|\nabla u|^2}g(JU, \nabla u) \right)\nonumber\\
&=&\frac{\nabla u}{|\nabla u|^2}\Hess \, u \left(U-\frac{g(U, \nabla
u)}{{|\nabla u|^2}}\nabla u, \nabla u\right)-\frac{J\nabla
u}{|\nabla u|^2}\Hess \, u \left(JU-\frac{g(JU, \nabla u)}{{|\nabla
u|^2}}\nabla u ,
\nabla u\right)\nonumber\\
&=&0, \ean where the last equality follows from \eqref{nabla} and
that $U-\frac{g(U, \nabla u)}{{|\nabla u|^2}}\nabla u\perp \nabla
u$.
\item $\gamma '=J\nabla u$: Using $J\nabla_X \nabla u=\nabla_{JX} \nabla u$ it reduces to the previous case.
\end{enumerate}

Now we have an orthogonal decomposition of the tangent bundle $TM=T_1\oplus T_1^{\perp}$ that is invariant under parallel transport.
By DeRham's decomposition theorem on a simply-connected manifold \cite[Page 187]{Kobayashi-Nomizu1963},  $M$ is a Riemannian product, and all the claims in the theorem follow.
\qed
\newline

\noindent {\em Proof of Corollary~\ref{coro}} Since the manifold $M$ is compact, by \cite{KK}, we can assume $\lambda >0$. Then, by \cite{Qian1997}, the universal cover $\tilde{M}$ is also compact. Now the result follows from Theorem~\ref{Kahler} and \ref{dim2}.
\qed

\end{document}